\documentclass[12pt]{amsart}
\usepackage{graphicx}

\usepackage{amsmath,graphics}
\usepackage{amsfonts,amssymb}
\usepackage[all]{xypic}

\theoremstyle{plain}
\newtheorem*{theorem*}{Theorem}
\newtheorem*{lemma*} {Lemma}
\newtheorem*{corollary*} {Corollary}
\newtheorem*{proposition*} {Proposition}
\newtheorem*{conjecture*}{Conjecture}

\newtheorem{theorem}{Theorem}[section]

{\catcode`@=11\global\let\c@figure=\c@theorem}

{\catcode`@=11\global\let\c@equation=\c@theorem}

\theoremstyle{remark}

\theoremstyle{definition}

\def\be{\begin{equation}}
\def\ee{\end{equation}}

\def\Q{\Bbb{Q}}

\def\Z{\Bbb{Z}}

\def\R{\Bbb{R}}

\def\N{\Bbb{N}}
\def\l{\lambda}

\def\part{\partial}

\def\bp{\begin{pmatrix}}
\def\sm{\smallsetminus}
\def\ep{\end{pmatrix}}
\def\bn{\begin{enumerate}}

\def\en{\end{enumerate}}
\def\ba{\begin{array}}
\def\ea{\end{array}}

\def\L{\Lambda}

\def\ti{\tilde}

\def\fr12{\frac{1}{2}}

\def\int{\mbox{Int}}

\def\ext{\mbox{Ext}}

\def\ol{\overline}
\def\sm{\setminus}

\newcommand{\Hom}{\operatorname{Hom}}

\begin{document}

\title{Realizations of Seifert matrices by hyperbolic knots}
\author{Stefan Friedl}

\address{Universit\'e du Qu\'ebec \`a Montr\'eal, Montr\'eal, Qu\`ebec}
\email{friedl@alumni.brandeis.edu}
% Instead of \subjclass[2000]{...},
% I put the above workaround for the old versions of amsart.
\def\subjclassname{\textup{2000} Mathematics Subject Classification}
\expandafter\let\csname subjclassname@1991\endcsname=\subjclassname
\expandafter\let\csname subjclassname@2000\endcsname=\subjclassname

\date{\today}
\begin{abstract}
Recently Kearton showed that any Seifert matrix of a knot is S--equivalent to the Seifert matrix of a prime knot. We show in this note that such a matrix is in fact S--equivalent to the Seifert matrix of a hyperbolic knot. This result follows from reinterpreting this problem in terms of Blanchfield pairings and by applying results of Kawauchi.
\end{abstract}

\maketitle

%===========================================
\section{Introduction}

We say that a square integral matrix  $A$ is of \emph{Seifert type} if $\det(A-A^t)=1$. Let $A$ be a square integral matrix, then for any column vector $v$ the matrices
\[ \bp A &0&0\\ v^t&0&0\\ 0&1&0\ep \mbox{ and } \bp A&v&0\\ 0&0&1\\ 0&0&0\ep\]
are called elementary enlargements of $A$. We also say that $A$ is an elementary reduction of any of its elementary enlargements. Two matrices are \emph{S--equivalent} if they can be connected by a chain of elementary enlargements, elementary reductions and unimodular congruences.

Let $K\subset S^3$ be a knot and $\emph{F}$ a Seifert surface. Given a basis for $H_1(F)$ we can then define the Seifert matrix $A$ of $K$. It is well--known that $A$ is of Seifert type. It is shown in \cite[Theorem~3.1]{Mu65} (cf. also \cite[Theorem~1]{Le70}) that the $S$--equivalence class of the Seifert matrix is a knot invariant.

It is well--known that any matrix of Seifert type is the Seifert matrix of a knot. In \cite{Ke04} Kearton showed that any matrix of Seifert type is S--equivalent to the Seifert matrix of a prime knot.

In this note we prove the following:

\begin{theorem}
Let $A$ be a matrix of Seifert type, then there exist infinitely many hyperbolic knots $K_i, i\in \N$ such that $A$ is S--equivalent to a Seifert matrix of $K_i$.
\end{theorem}

The proof relies on a reformulation of the S--equivalence class in terms of Blanchfield pairings and on realization results of Kawauchi.

Added in proof: This theorem also follows for links from combining Theorem \ref{thm:kawauchi} with
\cite[Theorem~A.1]{Ka94}.

%===========================================
\section{Proof of the theorem}

\subsection{S--equivalence and Blanchfield forms}
Given a knot $K\subset S^3$ we write $X(K)=S^3\sm \nu K$, the knot exterior. In the following we let $\L=\Z[t,t^{-1}]$ and $Q(\L)=\Q(t)$ the quotient field of $\L$. We view $\L=\Z[t,t^{-1}]$ with the involution $p\mapsto \ol{p}$ induced by $t\mapsto t^{-1}$.

Consider the following sequence of $\L$--homomorphisms
\[ \ba{rcl}
H_1(X(K);\L)&\xrightarrow{\cong}&H_1(X(K),\partial X(K);\L)\xrightarrow{\cong} H^2(X(K);\L)\xrightarrow{\cong}
\ext_\L^1(H_1(X(K);\L),\L)\\
&\xleftarrow{\cong}&\Hom(H_1(X(K);\L),Q(\L)/\L).\ea \]
Here the first map comes from the long exact sequence of the pair $(X(K),\partial X(K))$, and is easily seen to be an isomorphism. The second homomorphism is Poincar\'e duality, the third homomorphism comes from the universal coefficient spectral sequence (and is an isomorphism by \cite[Proposition~3.2]{Le77}) and finally the last homomorphism comes from the long exact Ext--sequence corresponding to the short exact sequence of coefficients
\[ 0\to \L\to Q(\L)\to Q(\L)/\L\to 0.\]
This sequence of homomorphisms defines the Blanchfield pairing
\[ \l(K): H_1(X(K);\L)\times H_1(X(K);\L)\to Q(\L)/\L.\]
This pairing is non--singular and $\L$--hermitian. Furthermore if $A$ is a Seifert matrix for $K$ of size $k\times k$, then the Blanchfield pairing is isometric to the pairing
\[ \ba{rcl}
\L^k/(At-A^t)\L^k\times \L^k/(At-A^t)\L^k&\to&Q(\L)/\L\\
(v,w)&\mapsto&\ol{v}^t(t-1)(At-A^t)^{-1}w.\ea \]
In particular the (S--equivalence class of a) Seifert matrix determines the Blanchfield pairing of a knot. By \cite{Tr73} (and also by comparing \cite{Ke75} with \cite{Le70}) the converse holds as well, more precisely, the following theorem holds true.

\begin{theorem} \label{thm:trotter}
Let $K_1,K_2\subset S^3$ be knot. Then $K_1$ and $K_2$ have S--equivalent Seifert matrices if and only if the Blanchfield pairings $\l(K_1)$ and $\l(K_2)$ are isometric.
\end{theorem}

\subsection{Kawauchi's realization results}
Before we continue we recall that the derived series $G^{(n)},n\in \N$ of a group $G$ is defined inductively by $G^{(0)}=G$ and $G^{(n+1)}=[G^{(n)},G^{(n)}]$, the commutator of $G^{(n)}$. We recall the following hyperbolic realization result by Kawauchi.

\begin{theorem}\label{thm:kawauchi}
Let $L\subset S^3$ be any link, then for any $V\in \R$ there exists a hyperbolic link
$\ti{L}\subset S^3$ together with a map $f:(S^3,\ti{L})\to (S^3,L)$ such that the following hold:
\bn
\item $\mbox{Vol}(S^3\sm \ti{L})>V$,
\item the induced map $\pi_1(S^3\sm \ti{L})/\pi_1(S^3\sm \ti{L})^{(n)}\to \pi_1(S^3\sm L)/\pi_1(S^3\sm L)^{(n)}$ is an
isomorphism for any $n$.
\en
\end{theorem}

The theorem follows from the theory of almost identical imitations of Kawauchi. More precisely,
the theorem follows from combining \cite[Theorem~1.1]{Ka89b} with
\cite[Properties~I~and~V,~p.~450]{Ka89a} (cf. also \cite{Ka89c}).

\subsection{Conclusion of the proof of the theorem}
Let $K\subset S^3$ be a knot and $V\in \R$. Let $\ti{K}$ be as in Theorem \ref{thm:kawauchi}. Since we can choose $V$ arbitrarily large it follows from Theorem \ref{thm:trotter} that it is enough to show that the Blanchfield pairings $\l(K)$ and $\l(\ti{K})$ are isometric.

First note that by Theorem \ref{thm:kawauchi} (2), applied to $n=1$, we have a commutative diagram
\[ \xymatrix{ \pi_1(X(\ti{K}))\ar[rr]^{f_*}\ar[dr] && \pi_1(X(K))\ar[dl]\\
&\Z.&}\]
In particular we get induced maps $H_i(X(\ti{K});\L)\to H_i(X({K});\L)$.
Write $X=X(K)$ and $\ti{X}=X(\ti{K})$. We then get the following commutative diagram
\[
\ba{ccccccccccccccccccccccccc}
H_1(\ti{X};\L)&\hspace{-0.3cm}\to&\hspace{-0.3cm}H_1(\ti{X},\partial \ti{X};\L)&\hspace{-0.3cm}\to&\hspace{-0.3cm}H^2(\ti{X};\L)&\hspace{-0.3cm}\to&\hspace{-0.3cm}\ext^1_{\L}(H_1(\ti{X};\L),\L)
&\hspace{-0.3cm}\xleftarrow{\cong}&\hspace{-0.3cm}\Hom(H_1(\ti{X};\L),Q(\L)/\L)\\[0.2cm]
\downarrow&\hspace{-0.3cm}&\hspace{-0.3cm}\downarrow&\hspace{-0.3cm}&\hspace{-0.3cm}\downarrow&\hspace{-0.3cm}&\hspace{-0.3cm}\downarrow&\hspace{-0.3cm}&\hspace{-0.3cm}\downarrow\\[0.2cm]
H_1({X};\L)&\hspace{-0.3cm}\to&\hspace{-0.3cm}H_1({X},\partial {X};\L)&\hspace{-0.3cm}\to&\hspace{-0.3cm}H^2({X};\L)&\hspace{-0.3cm}\to&\hspace{-0.3cm}\ext^1_{\L}(H_1({X};\L),\L)
&\hspace{-0.3cm}\xleftarrow{\cong}&\hspace{-0.3cm}\Hom(H_1({X};\L),Q(\L)/\L).\ea \]
This means that we get a commutative diagram
\[ \ba{ccccc}
H_1(X(\ti{K});\L)&\times &H_1(X(\ti{K});\L)&\to &Q(\L)/\L\\[0.1cm]
\downarrow&&\downarrow&&\downarrow =\\[0.1cm]
H_1(X({K});\L)&\times &H_1(X({K});\L)&\to &Q(\L)/\L. \ea \]
But it follows from Theorem \ref{thm:kawauchi} (2), applied to $n=2$, that the induced map $H_1(X(\ti{K});\L)\to H_1(X(K);\L)$ is an isomorphism of $\L$--modules. In particular $\l(\ti{K})$ is isometric to $\l(K)$.

\end{document}